\documentclass[12pt]{article}

\usepackage{a4wide}
\usepackage{amsmath}
\usepackage{enumerate}

\usepackage{amsmath,amsthm}
\usepackage{amsfonts}
\usepackage{amssymb}
\usepackage{longtable}
\usepackage[matrix,arrow,curve]{xy}

\newtheorem{theorem}{Theorem}
\newtheorem{cor}{Corollary}
\newtheorem{prop}{Proposition}

\newtheorem{lem}{Lemma}

\newtheorem{ex}{Example}{\bf}{\rm}

\def\Real{\mathbb{R}}

\def\D{\mathcal{D}}
\def\C{\mathcal{C}}

\def\F{\mathcal{F}}

\def\GLn{\text{\rm GL}(n,\Real)}
\def\On{\text{\rm O}(n)}
\def\GL1{\text{\rm GL}(1,\Real)}
\def\O1{\text{\rm O}(1)}
\def\GVL{\text{\rm GVL}}
\def\CL{\text{\rm CL}}
\def\Diff{\text{\rm Diff}}

\def\pt{\mathop\text{\rm pt}\nolimits}

\usepackage{stackrel}

\newcommand{\be}{\begin{equation}}
\newcommand{\ee}{\end{equation}}

\def\Im{\text{\rm Im}}

\usepackage{authblk}

\title{Losik classes for codimension-one foliations}

\author[1]{Yaroslav V. Bazaikin}
\author[2]{Anton S. Galaev}
\affil[1]{Sobolev Institute of Mathematics Novosibirsk, Russia and University of Hradec Kr\'alov\'e, Faculty of Science, Rokitansk\'eho 62, 500~03 Hradec Kr\'alov\'e,  Czech Republic
	E-mail: bazaikin(at)math.nsc.ru}
\affil[2]{University of Hradec Kr\'alov\'e, Faculty of Science, Rokitansk\'eho 62, 500~03 Hradec Kr\'alov\'e,  Czech
	Republic 
	E-mail: anton.galaev(at)uhk.cz}

\begin{document}

\maketitle

\begin{abstract}

Following Losik's  approach to Gelfand’s formal geometry, certain
characteristic classes  for codimension-one foliations coming from
the Gelfand-Fuchs cohomology are considered. Sufficient conditions for
non-triviality in terms of  dynamical properties of generators
of the holonomy groups are found. The non-triviality for the Reeb
foliations is shown; this is in contrast with some classical theorems
on the Godbillon-Vey class, e.g, the Mizutani-Morita-Tsuboi Theorem
about triviality of the Godbillon-Vey class of foliations almost
without holonomy is not true for the classes under consideration. It
is shown that the considered classes are trivial for a large class
of foliations without holonomy. The question of triviality is
related to ergodic theory of dynamical systems on the circle and to
the problem of smooth conjugacy of local diffeomorphisms. Certain
classes are obstructions for the existence of  transverse affine and
projective connections.
\vskip0.5cm


\vskip0.5cm


\tableofcontents

\end{abstract}

\section*{Introduction}

The Godbillon-Vey class \cite{GV} is an important invariant of a
codimension-one foliation on a smooth manifold $M$. It can be
explicitly given as the cohomology class $[\eta\wedge d\eta]\in
H^3(M)$ for a certain 1-form $\eta$ on $M$. Despite the simple
description, its geometrical meaning remains a mystery. There are
many works (e.g., \cite{Ref12,Ref26,Ref28,Ref82,Ref84,Ref106})
studying the relation of (non-)triviality of the Godbillon-Vey class
to geometry and dynamics of the foliation, see the surveys
\cite{Hurder00,Hurder16} and the book \cite{FolII} for the complete
list of references. It is natural to ask for stronger invariants
that could provide more information about  geometry of the
foliation, e.g., in \cite{Ref26,Hurder00}, it was suggested to study
the Vey class $[d\eta]$ introduced by Duminy and taking values in
certain cohomology associated to the foliation.

Losik~\cite{L90,L94,L15} introduced the notion of a generalized
atlas on the leaf space $M/\F$ of a foliation of codimension $n$
showing that there is a rich smooth structure on $M/\F$. In
particular, this approach allowed him to define the characteristic
classes as elements  of infinite order frame bundles
 over $M/\F$ coming from the generators of the
Gelfand-Fuchs cohomology of the Lie algebra $W_n$ of formal
 vector fields
on $\Real^n$. If the codimension $n$ of the foliation $\F$ is one, then the 
 generators of the cohomology $H^3(W_1,\O1)\cong H^2(W_1,\GL1)\cong\Real$ define  characteristic classes taking values in $H^3(S(M/\F)/\O1)$ and $H^2(S(M/\F)/\GL1)$, where $S(M/\F)$ is the infinite order frame bundle over the leaf space $M/\F$. We call these classes, respectively, the Godbillon-Vey-Losik class and the first Chern-Losik class
 (Losik called the last class the first Chern class). These classes
 may be interpreted as (secondary and primary) characteristic
 classes of the $\GL1$-bundle $S(M/\F)\to S(M/\F)/\GL1$.
 The classes may be mapped to the cohomology
of the manifold $M$. The image of the Godbillon-Vey-Losik class in
$H^3(M)$ is the classical Godbillon-Vey class. The first Chern-Losik
class always maps to zero in $H^2(M)$. One may consider the values
of these classes also in $H^3(S_2(M/\F)/\O1)$ and
$H^2(S_2(M/\F)/\GL1)$, where $S_2(M/\F)$ is the second order frame bundle over $M/\F$. Losik \cite{L15} proved that the first Chern class for the Reeb foliations
 on the three-dimensional sphere is non-trivial in $H^2(S_2(M/\F)/\GL1)$.
We define the Godbillon-Vey-Losik class with values in
$H^3(S(M/\F)/(\O1\times\mathbb{Z}))$. This class may be interpreted
as a secondary characteristic class of the $\mathbb{S}^1$-bundle
$S(M/\F)/(\O1\times\mathbb{Z})\to S(M/\F)/\GL1$.

We also define the Duminy-Losik class, which may be projected to the
Vey class. We prove

 {\bf Theorem.} {\it Let $\F$ be a codimension-one foliation on a manifold $M$. Then the Duminy-Losik class of $\F$
 is trivial if and only if $\F$ admits a transverse affine connection. }

We define also a characteristic class tightly related to the
existence of a transverse projective connection.
We find a sufficient condition for non-triviality of the Godbillon-Vey-Losik class with values in
$H^3(S(M/\F)/(\O1\times\mathbb{Z}))$
 and the first Chern-Losik class. This condition appeared to be deeply related to dynamical properties of the local diffeomorphisms generating holonomy of the foliation.

 {\bf Theorem.} {\it Let $\F$ be a codimension-one foliation on a manifold $M$. If the holonomy group of some leaf of $\F$ contains the germ of a local diffeomorphism
generating an infinite group and with the first derivative equal to
one at that point, then both the Godbillon-Vey-Losik class with
values in $H^3(S(M/\F)/(\O1\times\mathbb{Z}))$
 and the first Chern-Losik class of $\F$ are non-trivial. }

 The crucial fact that we use is the existence of a one-parameter group of local $C^1$-diffeomorphisms extending
 the iterations $\varphi^n$. In particular, we obtain
 
 {\bf Theorem.} {\it
The Godbillon-Vey-Losik class with values in
$H^3(S(M/\F)/(\O1\times\mathbb{Z}))$, and the first Chern-Losik class
for the Reeb foliations on the three-dimensional sphere are
non-trivial.}

This implies that the classical Mizutani-Morita-Tsuboi Theorem
\cite{Ref82} about triviality of the Godbillon-Vey class of
foliations almost without holonomy is not true for the classes we
consider. Recall that a foliation is called almost without holonomy
if the only leaves with non-trivial holonomy are compact. We see
that the Godbillon-Vey-Losik and  first Chern-Losik classes are more
sensitive to the holonomy than the
 Godbillon-Vey class. On the other hand, these classes vanish if the space of leaves is a one-dimensional orbifold, i.e., the classes do not
 contain the complete information about the holonomy.
 The case of the Reeb foliations shows that the classes under
 consideration does not satisfy the  theorem by Duminy and Sergiescu
\cite{Ref28} stating that a foliation without resilient leaves has
trivial Godbillon-Vey class. Similarly, the conjecture by Moussu
and Pelletier, and Sullivan that is  proved by Hurder
\cite{Hurder86} and stating that foliations with non-exponential
growth have trivial Godbillon-Vey class, is not true for the classes
that we consider.
In addition, unlike the classical Godbillon-Vey class, the classes under consideration are not invariant by cobordism.
Finally we study the question, whether the Morita-Tsuboi Theorem
\cite{Ref84}, stating that the Godbillon-Vey class is trivial for a
codimension-one foliation without holonomy on a compact manifold,
holds true for the classes under consideration. It is sufficient to
consider a foliation such that its leaf space is the factor of
 $\Real$ by a free Abelian group with generators $\varphi_0,\dots,\varphi_p\in\Diff_+(\Real)$ that have no fixed points.
 Consider a pair of diffeomorphisms
  $\varphi_i$ and $\varphi_j$. Factorizing $\Real$ by the group generated by $\varphi_i$, we obtain a circle; the diffeomorphism $\varphi_j$ acts on this circle
  and we may define its rotation number $\rho(\varphi_i, \varphi_j)$, which is irrational.

 {\bf Theorem.} {\it Let $\F$ be a codimension-one foliation without holonomy on a compact manifold $M$.
If for some $i$ and $j$ with $i\neq j$, the rotation number
$\rho(\varphi_i, \varphi_j)$ is Diophantine, then the Duminy-Losik
class, the Godbillon-Vey-Losik class and the first Chern-Losik class
of $\F$ are trivial.}

Thus we were able to generalize the Morita-Tsuboi Theorem almost for all foliations
 and the only case which remains unclear is the case of all numbers $\rho(\varphi_i, \varphi_j)$, $i\neq j$, being Liouville.
  This question is related to the problem of $C^\infty$-conjugation of the circle diffeomorphisms to rotations.

In the last section we give a picture comparing the triviality
conditions for the classes under consideration. Finally let us note that in the recent work \cite{BGG} we constructed examples of the Reeb foliations both with trivial and non-trivial Godbillon-Vey-Losik class with values in $H^3(S(M/\F)/\O1)$.

\section{Definition of characteristic classes following Losik}\label{sec1}

Here we review Losik's approach to the leaf spaces of foliations
\cite{L90,L94,L15} that allowed him to define new characteristic
classes of foliations. Let $\D_n$ be the category whose objects
are open subsets of $\Real^n$, and morphisms are \'etale (i.e.,
regular) maps. The dimension $n$ may be infinite; in this case we
use the definitions from \cite{BR} of manifolds with the model
space $\Real^\infty$. In that paper we use only the facts that the space 
$\Real^\infty$ is the limit of the projective sequence $0\leftarrow \Real\leftarrow\Real^2\leftarrow\cdots$, and the differential forms on 
$\Real^\infty$ may be locally represented as $\pi_k^*\omega$, where $k$ is a positive integer, $\pi_k:\Real^\infty\to\Real^k$ is the projection, and $\omega$ is a differential form on $\Real^k$.  

Let us recall the definition of a $\D_n$-space. Let $X$ be a set. A
$\D_n$-chart on $X$ is a pair $(U,k)$, where $U\subset \Real^n$ is an open subset, and
$k:U\to X$ is an arbitrary map. For two charts $k_i:U_i\to X$, a
morphism of charts is an \'etale map $m:U_1\to U_2$ such that
$k_2\circ m=k_1$. Let $\Phi$ be a set of charts and let $\C_\Phi$ be
the category whose objects are elements of $\Phi$ and morphisms are
all possible morphisms of the charts. The set $\Phi$ is called a
$\D_n$-atlas on $X$ if $X=\varinjlim J$, where $J:C_\Phi\to{\rm
Sets}$ is the obvious functor. A $\D_n$-space is a set $X$ with a maximal $\D_n$-atlas $\Phi$. A maximal atlas may be characterized by the following two conditions: 1) if $(V,k)\in\Phi$ and $m:U\to V$ is an \'etale map, then $(U,k\circ m)\in\Phi$; 2) if $U\subset\Real^n$ is an open subset, $\{U_i\}$ is an open cover of $U$, and $k:U\to X$ is a map such that $(U_i,k|_{U_i})\in\Phi$ for all $i$, then $(U,k)\in\Phi$.

If $\F$ is a foliation of codimension $n$ on a smooth manifold $M$,
then the leaf space $M/\F$ is a $\D_n$-space. The maximal
$\D_n$-atlas on $M/\F$ consists of the projections $U\to M/\F$,
where $U$ is a transversal which is the embedded to $M$ open subset
of $\Real^n$. These transversals may be obtained from a foliation
atlas on $M$.
Generally $\D_n$-spaces are orbit spaces of  pseudogroups of local
diffeomorphisms of smooth manifolds. Considering the space $\Real^n$
and the pseudogroup of all local diffeomorphisms of open subsets of
$\Real^n$, we see that the point $\pt$ is a $\D_n$-space. The atlas
of $\pt$ consists of all pairs $(U,k)$, where $U\subset\Real^n$ is an
open subset and $k:U\to \pt$ is the unique map. It is important to
note that $\pt$ is the terminal object in the category of
$\D_n$-spaces.

Each (co)functor from the category   $\D_n$ to the category of sets
may be extended to a (co)functor from the category of $\D_n$-spaces.
In this way one obtains, e.g., the de~Rham complex $\Omega^*(X)$ of
a $\D_n$-space $X$, which defines the de~Rham cohomology $H^*(X)$ of
$X$. E.g., if $X=M/\F$, then $\Omega^*(X)$ coincides with the
complex of basic forms.

Consider the functor $S$ assigning to each open subset
$U\subset\Real^n$ the space of frames of infinite order, i.e., the
space of jets at $0\in\Real^n$ of regular maps from $\Real^n$ to
$U$. Then for each $\D_n$-space $X$, we obtain the space $S(X)$ of frames
of infinite order. This space is a $\D_\infty$-space.
Similarly, let $S'(U)=S(U)/\GLn$ and $S''(U)=S(U)/\On$. We obtain
the $\D_\infty$-spaces $S'(X)$ and~$S''(X)$.

Let us consider the point $\pt$ as a $\D_n$-space. Each of the
spaces $S(\pt)$, $S'(\pt)$, and $S''(\pt)$ consists of a single
point; on the other hand, the complexes $\Omega^*(S(\pt))$,
$\Omega^*(S'(\pt))$, and $\Omega^*(S''(\pt))$, are naturally
isomorphic to the complexes $C^*(W_n)$, $C^*(W_n,\GLn)$, and
$C^*(W_n,\On)$, respectively, \cite{L94}. Here $W_n$ is the Lie algebra of
formal vector fields on $\Real^n$, and continuous cochains   with values in $\Real$ are considered. The just mentioned isomorphisms may be obtained using the canonical Gelfand-Kazhdan form with values in $W_n$ which is defined in the following way \cite{BR}. 
Let $U\subset\Real^n$ be an open subset. 
Let $\tau$ be a tangent
vector at $s\in S(U)$ and  let $s(u)$ be a curve on $S(U)$ such that
$\tau=\frac{ds}{du}(0)$. One can represent $s(u)$ by a smooth family
$k_u$ of germs at $0$ of regular at $0\in\Real^n$ maps $\Real^n\to U$,
i.e., $s(u)=j^\infty_0k_u$. Then  put
$$
\theta(\tau)=-j_0^\infty\frac{d}{du}(k_0^{-1}\circ k_u)(0).
$$
The form $\theta$ is canonical, i.e., it is invariant under the infinite order frame bundle maps induced by the diffeomorphisms of open subsets of $\Real^n$. Consequently, $\theta$ is a well-defined 1-form on $S(\pt)$ with values in $W_n$.
Now, each $c\in C^k(W_n)$ defines $\theta_c\in
\Omega^k(S(\pt))$ by the formula
\begin{equation}\label{mapforms}\theta_c(V_1,\dots, V_k)=c(\theta(V_1),\dots, \theta(V_k)),\end{equation}
where $V_1,\dots, V_k$ are vector fields on $S(\pt)$. The form $\theta$ satisfies the Maurer-Cartan equation, this implies that the map $c\mapsto \theta_c$ defines the homomorphism $C^*(W_n)\to \Omega^*(S(\pt))$ of the complexes. This homomorphism induces homomorphisms of the relative complexes.
Let now $X$ be a $\D_n$-space. The unique morphism
$$p_X:X\to\pt$$ of $\D_n$-spaces induces the characteristic homomorphisms
\begin{align}\label{chi}
\chi :H^*(W_n) &\cong H^*(S(\pt))\to H^*(S(X)), \\ \label{chi'}\chi' :H^*(W_n,\GLn)&\cong H^*(S'(\pt))\to
H^*(S'(X)), \\ \label{chi''}\chi'' :H^*(W_n,\On)&\cong H^*(S''(\pt))\to
H^*(S''(X)).\end{align} The images of the generators of the Gelfand-Fuchs cohomology under
these maps give Losik's characteristic classes.

Assume that $X=M/\F$.   Let $P(\F)$, $S(\F)$, $S'(\F)$, $S''(\F)$ be, respectively,  the frame bundle of the normal bundle of $\F$ over $M$, the infinite order frame bundle of the normal bundle of $\F$ over $M$, and its relative versions. The infinite dimensional manifolds $P(\F)$ and $S(\F)$ are homotopy equivalent. Both $S'(\F)$ and $S''(\F)$ are homotopy equivalent to $M$. There is a bundle map $S(\F)\to S(M/F)$ covering the projection $M\to M/\F$. This gives the composition of the homomorphisms
\begin{align}\label{chiA} & H^*(W_n)\to H^*(S(X)) \to H^*(S(\F))\cong      H^*(P(\F)),\\
\label{chi'A} & H^*(W_n,\GLn)\to H^*(S'(X))\to H^*(S'(\F))  \cong H^*(M),\\
\label{chi''A} & H^*(W_n,\On)\to H^*(S''(X))\to H^*(S''(\F)) \cong H^*(M).\end{align}
The  composition \eqref{chi''A} delivers the characteristic classes
of the foliation $\F$ \cite{L94}. The composition \eqref{chi'A} does not give
any new information, since it coincides with the composition
\begin{equation}\label{GLnOn} H^*(W_n,\GLn)\to H^*(W_n,\On)\to H^*(S''(X))\to H^*(M).\end{equation}
 If $P(\F)$ admits a section $s:M\to P(\F)$, then \eqref{chiA}
 together with the map
$s^*:H^*(P(\F))\to H^*(M)$ give the secondary characteristic classes of $\F$.

\section{Gelfand-Fuchs cohomology and canonical forms}

In this section we collect some known facts in order to describe explicitly the above constructions in the case $n=1$.

It is known \cite{Fuchs} that the complex
$C^*(W_1)$ is generated by the 1-forms $c_{r}$, $r=0,1,2\dots$,
where $$c_r(\xi)=\frac{d^r\xi^1}{(dt)^r}(0),\quad
\xi=\xi^1\frac{d}{dt}\in W_1.$$ 
The differential is given by 
\begin{equation}\label{dcr} dc_r=\sum_{k=0}^r
\left(\begin{smallmatrix}r\\
k\end{smallmatrix}\right)c_{r-k+1}\wedge c_k=c_{r+1}\wedge
c_0+(r-1)c_r\wedge c_1+\sum_{k=2}^{r-1}
\left(\begin{smallmatrix}r\\
k\end{smallmatrix}\right)c_{r-k+1}\wedge c_k.\end{equation}
In particular, $dc_1=c_2\wedge c_0$.
The complex $C^*(W_1,\O1)$ is spanned by the forms $1$ and $c_{i_1}\wedge\cdots\wedge c_{i_k}$ with odd $i_1+\cdots+i_k$,  while  the complex $C^*(W_1,\GL1)$ is spanned by $1$ and $c_2\wedge c_0$.
The only non-trivial cohomology groups in the positive degree are  $$
H^3(W_1)=H^3(W_1,\O1)=\Real[c_0\wedge c_1\wedge c_2],\quad
H^2(W_1,\GL1)=\Real [c_2\wedge c_0].$$

Let $U,W\subset\Real$ be  open subsets such that $0\in W$. Denote by $t$ the coordinate on $W$.
Let $f:W\to U$ be a regular map at the point $0$. The corresponding element of $S(U)$ has the
coordinates \begin{equation} \label{coord1} z_p =\frac{d^p
	f}{(dt)^p}(0),\quad p=0,1,2,\dots\end{equation} Suppose that we have
another open subset $V\subset\Real$ and a regular map $h:U\to V$.
Denote by $\alpha_p$ the coordinates on $S(V)$. Then we get the
induced map $\tilde h:S(U)\to S(V)$ given by
\begin{equation}\label{coordtrans1}\alpha_0=h(z_0),\quad \alpha_n=n!\sum_{k=1}^n\frac{h^{(k)}(z_0)}{k!}\sum_{i_1+\dots+
	i_k=n}\frac{z_{i_1}}{i_1!}\dots\frac{z_{i_k}}{i_k!},\quad n\geq
 1.\end{equation} It is convenient to consider the following
coordinates on $S(U)$:
\begin{equation} \label{coord2} y_0 =z_0,\quad y_1=z_1,\quad y_2=\frac{z_2}{z_1^2},\quad  y_3=\frac{z_3}{z_1^3},\dots \end{equation}

The action of $\GL1=\Real^*$ on $S(U)$ is given by
\begin{equation}\label{Raction}(z_0,z_1,z_2,\dots)\mapsto (z_0,\lambda z_1,\lambda^2 z_2,\dots),\quad \lambda\in \Real^*.\end{equation}
With respect to the coordinates \eqref{coord2} it takes the form
\begin{equation}\label{Ractiony}(y_0,y_1,y_2,y_3,\dots)\mapsto (y_0,\lambda y_1,y_2, y_3,\dots),\quad \lambda\in \Real^*.\end{equation}

The functions
\begin{equation} \label{coord3} x_0 =y_0,\quad x_1=\ln|y_1|,\quad x_2=y_2,\quad  x_3=y_3,\dots \end{equation}
may be considered as the coordinates on $S''(U)=S(U)/\O1$.  The
coordinate transformation  takes the form
\begin{equation}\label{coordtrans2}\begin{aligned}\beta_0&=h(x_0),\quad \beta_1=x_1+\ln|h'(x_0)|,\\ \beta_n&=n!\sum_{k=1}^{n-1}\frac{1}{k!}\frac{h^{(k)}(x_0)}{(h'(x_0))^n}\sum_{i_1+\dots+
	i_k=n}
\frac{x_{i_1}}{i_1!}\dots\frac{x_{i_k}}{i_k!}+\frac{h^{(n)}(x_0)}{(h'(x_0))^n},\quad
n\geq 2,\end{aligned}\end{equation} where in the last equality it is
assumed that $x_1=1$. In particular,
\begin{equation}\label{coordtrans2A}
\beta_2=\frac{x_2}{h'(x_0)}+\frac{h''(x_0)}{(h'(x_0))^2},\quad
\beta_3=\frac{x_3}{(h'(x_0))^2}+\frac{3h''(x_0)x_2}{(h'(x_0))^3}+\frac{h'''(x_0)}{(h'(x_0))^3}.\end{equation}
 The functions $x_0,\quad x_2,\quad  x_3,\dots $ are coordinates on $S'(U)$, and the  projection
$$p:S''(U)\to S'(U)$$
is given by $$(x_0,x_1,x_2,x_3,\dots)\mapsto (x_0,x_2,x_3,\dots).$$

Consider the Gelfand-Kazhdan  form $\theta$ on $S(\pt)$. The equality $$\theta=\sum_{k\geq 0}\frac{1}{k!}\theta_k t^k\frac{d}{dt}$$ defines on $S(\pt)$ canonical 1-forms $\theta_0,\theta_1,\theta_2,\dots$. Applying \eqref{mapforms}, we get 
$$ \theta_k=\theta_{c_k},\quad k=0,1,2,\dots .$$
 From the definition of $\theta$ it follows that with respect to the coordinates \eqref{coord2} it holds  
\begin{equation}\label{formytheta}
\theta_0=-\frac{dy_0}{y_1},\quad\theta_1=-\frac{dy_1}{y_1}+y_2dy_0,\quad\theta_2=-y_1dy_2+y_1\big(y_3-2y_2^2\big)dy_0.\end{equation}
The canonical forms over smooth manifolds were introduced by Kobayashi \cite{Kob61} in order to extend the theory of the affine connections to the higher order frame bundles. In particular, the formulas   \eqref{formytheta} may be easily deduced from \cite{Kob61}. For more information on canonical forms including the relation to the formal vector fields see, e.g., \cite[Sec. 3.3]{Morita}.

\section{Godbillon-Vey-Losik class and the first Chern-Losik class}\label{secdefdim1}

Let $\F$ be a codimension-one foliation on a manifold $M$.
Below in this section we will see that the image of the class $[c_0\wedge c_1\wedge c_2]\in H^3(W_1,\O1)$ in $H^3(M)$ under the composition \eqref{chi''A} coincides with the Godbillon-Vey class of the
foliation $\F$. By that reason, for   a  $\D_1$-space $X$, we call the image of the class $[c_0\wedge c_1\wedge c_2]$ in $H^3(S''(X))$ under the map \eqref{chi''} the
 {\bf Godbillon-Vey-Losik class} (GVL class) of $X$
 and denote it by $\GVL(X)$. We see that the GVL class is represented by the form $$gvl=\theta_0\wedge\theta_1\wedge\theta_2\in\Omega^3(S''(X)).$$ 
 Note that \begin{equation}\label{rel2}gvl=\theta_1\wedge d\theta_1.\end{equation}
  With respect  
to the coordinates \eqref{coord3} it holds
\begin{equation*} gvl=-dx_0\wedge dx_1\wedge dx_2.\end{equation*}

The composition \eqref{GLnOn} shows that the image of $[c_2\wedge c_0]\in H^2(W_1,\GL1)$ in $ H^2(M)$ is always trivial.
On the other hand, the image of $[c_2\wedge c_1]$ in $H^2(S'(X))$ under the map \eqref{chi'} may be non-trivial \cite{L90} (also in the case $X=M/\F$).
We call this image {\bf the first Chern-Losik class} (first CL class) and denote it by $\CL_1(X)$. The first CL class is represented by the form
$$cl_1=\theta_2\wedge\theta_0\in\Omega^2(S'(X)).$$ With respect  
to the coordinates \eqref{coord3} it holds
\begin{equation*} cl_1=dx_2\wedge  dx_0.\end{equation*}

The first CL class and the GVL class have the following geometric
interpretation. The projection $S(X)\to S'(X)$ may be considered as
a principal $\GL1$-bundle. The form $\theta_1$ is a connection form
on that bundle. The curvature of this connection is given by
$d\theta_1=\theta_2\wedge \theta_0$. This form is projectable to the
base $S'(X)$ and it defines the first CL class. We see that the
first CL class may be obtained from the Chern-Weil construction. The
GVL class belongs to the cohomology of the total space of the
bundle; it is represented by the product of the connection form and
the curvature form, i.e., this class is the secondary characteristic
class in the sense of the Chern-Simons construction.

Let us explicitly describe the relation of the GVL class to the
classical Godbillon-Vey class. Let $\F$ be a foliation of
codimension one on a smooth manifold $M$. Suppose that the foliation
$\F$ is defined by a non-vanishing $1$-form $\omega$. Then there is a 1-form $\eta$ such that 
\begin{equation}\label{dweta}
d\omega=\eta\wedge\omega,
\end{equation}
and the Godbillon-Vey class of $\F$ is defined as $[\eta\wedge d\eta]\in H^3(M).$
The form $\eta$ may be chosen in the following way.
Let $g$ be a
Riemannian metric on $M$. Let $X$ be the vector field on $M$
orthogonal to $\F$ and such that $\omega(X)=1$. The form
$\eta=-L_X\omega$ satisfies $d\omega=\eta\wedge \omega$.  Let $x\in M$. Let
$(v_1,\dots,v_{n-1},u)$ be foliated coordinates around the point
$x$, where $u$ is the transverse coordinate. Let $\gamma(t)$ be the
integral curve of the vector field $X$ with $\gamma(0)=x$.
Considering the jet at $t=0$ of the function $u(t)$, we obtain the
map
$$\sigma:M\to S(M/\F).$$ Let $f$ be a non-vanishing function such
that in the coordinate neighborhood it holds
$\omega=\frac{1}{f}du$. It is not hard to check that, with respect
to the coordinates \eqref{coord1}, the map $\sigma$ is given by
$$(v_1,\dots,v_{n-1},u)\mapsto (u,f,Xf,*),$$
where the star stands for the values of the coordinates
$z_3,z_4\dots$. Then it holds
\begin{equation}\label{sigmatheta0}\sigma^*\theta_0=-\omega,\end{equation}
\begin{equation}\label{sigmatheta1}\sigma^*\theta_1=\eta.\end{equation}

 Let $\sigma'$ and
$\sigma''$ be the compositions of $\sigma$ with the projections
$S(M/\F)\to S'(M/\F)$ and $S(M/\F)\to S''(M/\F)$, respectively.
Then we see that
 $$\sigma''^*(gvl)=\sigma''^*(\theta_1\wedge d\theta_1)=\eta\wedge d\eta,$$
 i.e., the map $\sigma''^*$ maps the GVL class to the Godbillon-Vey
 class.
 Next, $$\sigma'^*(cl_1)=\sigma'^*(\theta_2\wedge\theta_0)=\sigma^*(d\theta_1)=d\eta,$$
 which confirms the above statement that the first CL class is
 always trivial in $H^2(M)$.

Note that the characteristic classes under consideration belong to the cohomology of the second order frame bundles $S'_2(X)$ and  $S_2''(X)$, which are $\D_2$ and $\D_3$-spaces, respectively. 
The projections $$S''(X)\to S''_2(X),\quad S'(X)\to S'_2(X)$$
induce the maps $$H^*(S''_2(X))\to H^*(S''(X)),\quad H^*(S'_2(X))\to H^*(S'(X))$$
sending characteristic classes to characteristic classes.
A priori these maps must not be injective. This shows that
the GVL class and the first CL class, considered  as elements of $H^3(S_2''(X))$ and $H^2(S_2'(X))$ respectively, may contain more information. In \cite{L15}, Losik
proved that the class $[\theta_2\wedge \theta_0]$ is non-trivial in $H^2(S'_2(M/\F))$ for the Reeb foliations  on the
three-dimensional sphere. Below we will prove a stronger
statement about the non-triviality of this class in $H^2(S'(M/\F))$.

\section{Godbillon-Vey-Losik class with value in
$H^3(S''(X)/\mathbb{Z})$}

Let $X$ be a $\D_1$-space. Consider the space $S''(X)$. Consider the restriction of the action \eqref{Ractiony} to the subgroup $\left\{e^n|n\in\mathbb{Z}\right\}\subset\Real\setminus\{0\}=\GL1$. With respect to the 
coordinates \eqref{coord3}, we obtain the action of
the group $\mathbb{Z}$ on $S''(U)$: $$(x_0,x_1,x_2,x_3,\dots)\mapsto (x_0,x_1+n,x_2,x_3,\dots),\quad n\in\mathbb{Z}.$$ The
projection $S''(U)\to S''(U)/\mathbb{Z}$ defines a $\D_\infty$-chart
on $S''(U)/\mathbb{Z}$. The actions of the group $\mathbb{Z}$
commute with the morphisms of charts from $S''(X)$. This shows that
we get a well-defined $\D_\infty$-space $S''(X)/\mathbb{Z}$.

\begin{prop}\label{propisomHS/Z} It holds $H^*(S''(\pt)/\mathbb{Z})\cong H^*(W_1)$.
\end{prop}

{\bf Proof.}   Consider coordinates
\eqref{coord1} on $S(\pt)$ and the action of the group
$\mathbb{Z}$ on $S(\pt)$ defined by
$$(z_0,z_1,z_2,\dots)\mapsto (z_0,e^nz_1,e^{2n}z_2,\dots),\quad n\in\mathbb{Z}.$$
This action commutes with the action of $\O1$, and it corresponds
to the action of $\mathbb{Z}$ on $S''(\pt)$. Using the definition
of the form $\theta$, it is not hard to see that
$$1\cdot \theta_{c_r}=\theta_{e^{r-1}c_r},$$
where $1\cdot$ denotes the action of the element $1\in\mathbb{Z}$.
This shows that the complex $\Omega^*(S(\pt)/\mathbb{Z})$ is spanned
by the forms $$1,\quad \theta_{1},\quad \theta_{2}\wedge
\theta_{0}, \quad \theta_{0}\wedge\theta_{1}\wedge
\theta_{2}.$$ These forms are invariant under the action of the
group $\O1$, i.e., these forms span the complex
$\Omega^*(S''(\pt)/\mathbb{Z}).$ It holds $d\theta_1=\theta_2\wedge \theta_0$. This
proves the statement of the proposition. \qed

Let now $X$ be an arbitrary $\D_1$-space. The map $X\to\pt$ defines
the homomorphism
$$H^*(S''(\pt)/\mathbb{Z})\to H^*(S''(X)/\mathbb{Z}),$$
and we obtain {\bf the GVL class with values in} $H^*(S''(X)/\mathbb{Z})$,
$$[\theta_0\wedge \theta_1\wedge \theta_2]\in
H^3(S''(X)/\mathbb{Z}).$$

The geometric interpretation of the just defined GVL class is the following. Consider the $\mathbb{S}^1$-bundle
$$S''(X)/\mathbb{Z}\to S'(X).$$ The form $\theta_1$ is a connection
form on that bundle. The GVL class is a secondary characteristic
class of the bundle. Again, the first CL class is a characteristic
class of that bundle.

The projection $S''(X)\to S''(X)/\mathbb{Z}$ induces the
homomorphism
$$H^3(S''(X)/\mathbb{Z})\to H^3(S''(X)),$$
which  preserves the GVL classes. This implies

\begin{prop}\label{proptrivial0} Suppose that for a $\D_1$-space $X$,  the GVL class is trivial in
$H^3(S''(X)/\mathbb{Z})$ (resp., in $H^3(S_2''(X)/\mathbb{Z})$).
Then it is trivial in $H^3(S''(X))$ (resp., in $H^3(S_2''(X))$).
 \end{prop}

The results from below and \cite{BGG} show that for some Reeb foliations, the GVL class is trivial in  $H^3(S_2''(X))$ and $H^3(S''(X))$, while it is not trivial in 
$H^3(S''(X)/\mathbb{Z})$ and $H^3(S_2''(X)/\mathbb{Z})$.
The next proposition gives a comparison of the GVL and the first
CL classes.

\begin{prop}\label{proptrivial} Suppose that for a $\D_1$-space $X$,  the first CL class is trivial in $H^2(S'_2(X))$.
Then the first CL class is trivial in $H^2(S'(X))$;
 the GVL class is trivial in $H^3(S_2''(X)/\mathbb{Z})$, and consequently also in $H^3(S''(X)/\mathbb{Z})$.
 \end{prop}
{\bf Proof.} Suppose that the first CL class is trivial in
$H^2(S'_2(X))$, then there is a 1-form
$\gamma\in\Omega^1(S'_2(X))$ such that $d\gamma=cl_1$. The first
statement of the proposition is obvious. Consider the projection
$$\pi:S''_2(X)/\mathbb{Z}\to S'_2(X).$$
Using \eqref{rel2}, we get
\begin{multline*}d(  \theta_1 \wedge \pi^*\gamma)=
d\theta_1\wedge \pi^*\gamma-\theta_1\wedge d\pi^*\gamma \\ =
\pi^*cl_1\wedge \pi^*\gamma-\theta_1\wedge
d\theta_1=-\theta_1\wedge d\theta_1+\pi^*(cl_1 \wedge
\gamma)=-gvl,\end{multline*} since $cl_1 \wedge
\gamma\in\Omega^3(S_2'(X))=0$. This implies the proof of the
proposition. $\Box$

Let $X$ be  a $\D_1$-space. Consider the projection
$$\pi:S''(X)/\mathbb{Z}\to S'(X).$$ It can be viewed as an
$\mathbb{S}^1$-bundle.
 Let $U\subset\Real$ be an open subset.
For a $k$-form
$$\omega=\sum_{i_1,\dots, i_k\neq 1}\omega_{i_1\cdots
i_k}dx_{i_1}\wedge\cdots\wedge dx_{i_k}\in
\Omega^k(S''(U)/\mathbb{Z}),$$ let
$$\pi_*\omega=0,\quad \pi_*(dx_1\wedge \omega)=-\sum_{i_1,\dots, i_k\neq 1}\left(\int_0^1\omega_{i_1\cdots
i_k}dx_1\right)dx_{i_1}\wedge\cdots\wedge
dx_{i_k}\in\Omega^k(S'(U)).$$ For any $\D_1$-space $X$, we obtain
the map
$$\pi_*:\Omega^*(S''(X)/\mathbb{Z})\to \Omega^{*-1}(S'(X)).$$
This map is usually called the Gysin homomorphism or the integration
along the fiber; this map  commutes with the exterior derivative.
Note that
$$\pi_*(\theta_0\wedge \theta_1\wedge \theta_2)=\theta_2\wedge\theta_0.$$

\begin{prop}\label{proptrivial1} Let $X$ be  a $\D_1$-space. If the GVL class of $X$
is trivial in $H^3(S''(X)/\mathbb{Z})$ (resp., in $H^3(S_2''(X)/\mathbb{Z})$), then the first CL class of
$X$ is trivial in $H^2(S'(X))$ (resp., in $H^2(S_2'(X))$). \end{prop}

{\bf Proof.} Suppose that $\omega\in\Omega^2(S''(X)/\mathbb{Z})$
and
$$d\omega =\theta_0\wedge \theta_1\wedge \theta_2.$$
 Since
$\pi_*$ commutes with the exterior derivative, we get
$$d\pi_*\omega=\theta_2\wedge\theta_0.$$
\qed

\section{Duminy-Losik class}

Let us recall the definition of the Vey class introduced by Duminy.
Let $\F$ be a foliation of codimension one on a  manifold $M$.
Suppose that the foliation $\F$ is defined by a non-vanishing
$1$-form $\omega$. Consider the complex
$$A^{m}=\Omega^{m-1}(M)\wedge \omega$$ with the differential being the usual exterior derivative.
Denote by $H^*_\omega(M)$ the cohomology of this complex. Let $\eta$
be any 1-form such that $d\omega=\eta\wedge \omega$. Then the Vey
class is the class of the form $d\eta$ in $H^2_\omega(M)$.

Let now $X$ be a $\D_1$-space. The canonical form $\theta_0$
satisfies
$$d \theta_0=\theta_1\wedge \theta_0.$$ Consider the complex
$$C^m=\Omega^{m-1}(S(X))\wedge \theta_0$$ with the differential
being the usual exterior derivative. Denote by
$H^*_{\theta_0}(S(X))$ the cohomology of this complex.

\begin{prop} It holds
$$H^2_{\theta_0}(S(\pt))=\Real [\theta_2\wedge \theta_0],\quad H^3_{\theta_0}(S(\pt))=\Real [\theta_0\wedge
\theta_1\wedge\theta_2],\quad H^k_{\theta_0}(S(\pt))=0,\quad k\neq
2,3.$$\end{prop}

{\bf Proof.}  It is clear that the complex
$\Omega^{*-1}(S(\pt))\wedge \theta_0$ is isomorphic to the complex
$C^{*-1}(W_1)\wedge c_0$. 

For each $k\geq 2$ and $0\leq s\leq k$, consider the complex
$L^k_s$ consisting of the forms
$$\lambda_k\wedge\cdots\wedge\lambda_{k-s+1}\wedge\alpha\wedge c_0,$$
where $\alpha$ is a combination of $c_1,\dots, c_{k-s}$, and it
holds $$d\lambda_r=a_r\lambda_r\wedge c_1$$ for some specific
constants $a_r>0$. The equality \eqref{dcr} shows that $L^k_0$ is
a subcomplex of $C^{*-1}(W_1)\wedge c_0$. Moreover, the complex $C^{*-1}(W_1)\wedge c_0$ is the injective limit of the complexes~$L^k_0$.

For a given $L^k_s$, let $\tilde L^k_{s}$ denote the subcomplex
consisting of the above elements such that $\alpha$ is a combination
of $c_1,\dots, c_{k-s-1}$. It is clear that $\tilde L^k_{s}$ is just
$L^{k-1}_{s}$ up to the choice of the constants~$a_r$.

\begin{lem} It holds
$$H(L_0^k)=\Real[c_2\wedge c_0]\oplus\Real[c_0\wedge c_1\wedge
c_2],\quad k\geq 2,$$
$$H^*(L_s^k)=0,\quad k\geq 3,\,\, 1\leq s\leq k-1.$$\end{lem}

{\bf Proof.} We prove the statement of the lemma by induction
over $k$ and $s$. The computation of the cohomology $H^*(L^3_1)$
and $H^*(L^2_0)$ is direct. Suppose that $k_0$ is fixed and the
statement holds true for all $k\leq k_0$ and all $1\leq s  \leq
k_0-1$. Let us prove the statement for $k=k_0+1$. Let us compute
the cohomology of the complex $L^{k_0+1}_{k_0}$. This complex is
spanned by the forms
$$\alpha=\lambda_{k_0+1}\wedge\cdots\wedge \lambda_2\wedge c_0,\quad \beta=\lambda_{k_0+1}\wedge\cdots\wedge \lambda_2\wedge c_1\wedge
c_0.$$ It holds $d\alpha=\pm (a_{k_0+1}+\cdots +a_2)\beta$. This
shows that $H^*(L^{k_0+1}_{k_0})=0$.
Suppose that the statement holds true for all $s$ such that
$s_0\leq s\leq k_0$ for some fixed $s_0$, $2\leq s_0\leq k_0-1$.
Let us prove the statement for $s=s_0-1$. Consider the short exact
sequence $$0\rightarrow \tilde L^{k_0+1}_{s_0}\rightarrow
L^{k_0+1}_{s_0}\rightarrow L^{k_0+1}_{s_0}/\tilde
L^{k_0+1}_{s_0}\rightarrow 0.$$ According to the induction
hypothesis, $H^*(L^{k_0+1}_{s_0})=0$, and
$0=H^*(L^{k_0}_{s_0})=H^*(\tilde L^{k_0+1}_{s_0})$. Note that
$L^{k_0+1}_{s_0}/\tilde L^{k_0+1}_{s_0}$ is just
$L^{k_0+1}_{s_0-1}$ for a proper choice of the constants $a_r$.
From the above exact sequence it follows that
$H^*(L^{k_0+1}_{s_0-1})=0$. Similarly,  from the short exact
sequence
$$0\rightarrow \tilde L^{k_0+1}_{0}\rightarrow
L^{k_0+1}_{0}\rightarrow L^{k_0+1}_{0}/\tilde
L^{k_0+1}_{0}\rightarrow 0$$ and the induction hypothesis it
follows that $H^*(L^{k_0+1}_{0})\cong H^*(L^{k_0}_{0})$.  \qed

Now, each element of the complex  $C^{*-1}(W_1)\wedge c_0$ belongs
to $L^k_0$ for some $k$. This completes the proof of the
proposition. \qed

Let $X$ be a $\D_1$-space. The characteristic morphism $X\to{\pt}$
induces the map
$$H^*_{\theta_0}(S(\pt))\to H^*_{\theta_0}(S(X)).$$
We call the  class $[\theta_2\wedge \theta_0]\in
H^2_{\theta_0}(S(X))$ {\bf the Duminy-Losik class} (DL class)
of~$X$.

For an arbitrary $\D_1$-space $X$, we may consider also the class
$[\theta_0\wedge \theta_1 \wedge \theta_2]\in H^3_{\theta_0}(S(X))$.
Note that if the DL class is trivial, then this class is trivial as
well. Indeed, if
$$d(F\theta_0)=\theta_2\wedge\theta_0,$$ then $$d(F\theta_1\wedge
\theta_0)=-\theta_0\wedge \theta_1 \wedge \theta_2.$$ Similarly, for
a given foliation $\F$, one may consider the cohomology class of the
form $\eta\wedge d\eta$ in $H^3_\omega(M)$. The triviality of the
Vey class implies the triviality of that class. 

Let $X$ be a $\D_1$-space. We say that an affine connection on $X$ is defined if an affine connection is fixed on the domain of each local chart on $X$ and the morphisms of charts are affine mappings. If $X$ is the leaf space of a foliation, this definition coincides with the definition of a transverse affine connection.    An affine connection on $X$  puts in
correspondence to each coordinate $w$ on $X$ the Christoffel symbol, which is a function $T_w$ of the
variable $w$. For each morphism of charts $t=\varphi(w)$ it
holds
$$T_t(\varphi(w))+\frac{\varphi''(w)}{(\varphi'(w))^2}=\frac{T_w(w)}{\varphi'(w)}.$$

\begin{theorem}\label{ThDLclass} The Duminy-Losik class of $X$ is trivial if and
only if $X$ admits an affine connection.
\end{theorem}
{\bf Proof.} The DL class of $X$ is trivial if and only there exists
a function $f$ on $S(X)$ such that
$$d(f\theta_0)=\theta_2\wedge\theta_0.$$
Equivalently, for each coordinate $y_0$ on $X$, it holds
$$f=-y_1(y_2+F_{y_0}(y_0))$$ for some function $F_{y_0}$;
moreover, under the coordinate transform $t_0=\varphi(y_0)$, it
holds
$$F_{t_0}(\varphi(y_0))+\frac{\varphi''(y_0)}{(\varphi'(y_0))^2}=\frac{F_{y_0}(y_0)}{\varphi'(y_0)}$$
(the last equality is equivalent to the fact that $f$ is a function
on $S(X)$). Here we used coordinates~\eqref{coord2}. The functions
$F_{y_0}$ define an affine connection  on $X$. \qed

The existence of an affine connection depends only on the frame bundle of the second order. This implies

\begin{cor} The DL class is
trivial in $H^2_{\theta_0}(S(X))$ if and only if it is trivial in
$H^2_{\theta_0}(S_2(X))$.  \end{cor}

It is clear that if the DL class of a $\D_1$-space is trivial, then the first CL class is trivial  both in  $H^2(S''(X))$ and $H^2(S_2''(X))$.
From this and Propositions \ref{proptrivial0} and  \ref{proptrivial} we obtain

\begin{cor}
If the DL class of $X$ is trivial, then all versions of the GVL class and the first CL class are trivial.\end{cor}

 Let again $\F$ be a foliation of
codimension one on a  manifold $M$ defined by a non-vanishing
$1$-form $\omega$. We use the notation from Section
\ref{secdefdim1}. The equality \eqref{sigmatheta0} shows that the
map $\sigma^*$ induces the map
$$H^*_{\theta_0}(S(M/\F))\to H^*_\omega(M),$$
which maps the DL class onto the Vey class. In particular, if the DL
class is trivial, then  the Vey class is trivial as well.

\begin{cor} If the foliation $\F$ admits a transverse affine
connection, then its Vey class is trivial. \end{cor}

Note that the Vey class of the Reeb foliations is trivial, while the
DL class of the Reeb foliations is non-trivial (this follows from the
results of the next section).

Let $X$ be a $\D_1$-space. Similar to the cohomology
$H^*_{\theta_0}(S(X))$, let us consider the cohomology
$H^*_{\theta_1\wedge\theta_0}(S(X))$ of the complex
$\Omega^{*-2}(S(X))\wedge\theta_1\wedge\theta_0.$

We say that a projective connection $q$ on $X$ is given if, for
each $\D_1$-chart with the coordinate $w$, a function $q_w(w)$ is
defined such that for each morphism of charts $t=\varphi(w)$, it
holds \begin{equation}\label{projcon}
q_t(\varphi(w))\cdot(\varphi'(w))^2=q_w(w)-{\rm
S}(\varphi)(w),\end{equation} where
$${\rm S}(\varphi)(w)=\frac{\varphi'''(w)}{\varphi'(w)}-\frac{3}{2}\frac{(\varphi''(w))^2}{(\varphi'(w))^2}$$
is the Schwarzian derivative.

Consider the canonical form
$$\theta_3=-y_1^2dy_3+3y_2y_1^2dy_2+y_1^2(y_4+6y_2^3-6y_2y_3)dy_0.$$
It is not hard to check that the form
$\theta_3\wedge\theta_1\wedge\theta_0$ defines a non-trivial class
in $H^3_{\theta_1\wedge\theta_0}(S(\pt))$. This class defines a
characteristic class for each $\D_1$-space.

\begin{prop} Let $X$ be a $\D_1$-space. Then the class $[\theta_3\wedge\theta_1\wedge\theta_0]$ is
trivial in $H^3_{\theta_1\wedge\theta_0}(S(X))$ if and only if  for each coordinate $y_0$ on $X$ a function
 $q_{y_0}(y_0,y_1)$ on $S_1(U)$ is given such that
under the coordinate transformation $t_0=\psi(y_0)$ it holds
\begin{equation} \label{eqq}q_{t_0}(\varphi(y_0),\varphi'(y_0)y_1)(\varphi'(y_0))^2=q_{y_0}(y_0,y_1)-{\rm
S}(\varphi)(y_0).\end{equation} Consequently, if $X$ admits a projective connection, then the class $[\theta_3\wedge\theta_1\wedge\theta_0]$ is
trivial in $H^3_{\theta_1\wedge\theta_0}(S(X))$.
\end{prop}

{\bf Proof.} The triviality of the class under  consideration is equivalent to the existence of a function $F$ on $S(X)$ such that
$$d(F \theta_1\wedge\theta_0)=\theta_3\wedge\theta_1\wedge\theta_0.$$
With respect to local coordinates \eqref{coord2}, the last equality is equivalent to
$$dF\wedge d y_1\wedge dy_0=\theta_3\wedge dy_1\wedge d y_0.$$
Equivalently,
$$F=-y_1^2y_3+\frac{3}{2}y_1^2y_2^2+G_{y_0}(y_1,y_0)$$ for some function $G_{y_0}(y_1,y_0)$ associated to the coordinate $y_0$ on $X$.
Since $F$ is a function on $S(X)$, we get that under any coordinate transformation $t_0=\psi(y_0)$, it holds
$$G_{t_0}(\varphi(y_0),\varphi'(y_0)y_1)=G_{y_0}(y_0,y_1)-y_1^2{\rm
S}(\varphi)(y_0).$$ Equivalently, the functions $$q_{y_0}(y_0,y_1)=\frac{G_{y_0}(y_0,y_1)}{y_1^2}$$ satisfy \eqref{eqq}. This proves the proposition. \qed

\section{Non-triviality for the Reeb foliations and consequences}

In \cite{L15}, Losik proved that the class
$[\theta_2\wedge\theta_0]$ is  non-trivial in
$H^2(S_2'(\mathbb{S}^3/\F))$ for the Reeb foliations on the
3-dimensional sphere. The next theorem makes this statement
stronger. The statement of the theorem follows also from more
general Theorem \ref{Nontrivil1A} proven below.

\begin{theorem}\label{ThCLReeb} The  first Chern-Losik class is non-trivial in $H^2(S'(\mathbb{S}^3/\F))$ for the Reeb foliations on the 3-dimensional sphere. \end{theorem}

{\bf Proof.} Let us recall the description of the Reeb foliations.
Consider a smooth function $f(x)$  on the interval $(-1,1)$
satisfying the following conditions:
\begin{align}\label{f}
&f(0)=0,\,f(x)\ge 0,\,f(-x)=f(x),\\ \label{f1}\lim_{x\to\pm
1}&\frac{d^pf}{(dx)^p}(x)=\infty,\,\lim_{x\to\pm
1}\frac{d^p}{(dx)^p}\frac1{f'(x)}=0,\,\mbox{for $p=0,1,\dots$}
\end{align}
In \cite{L15}, it is shown that the function $f$ has the
following properties:
\begin{equation}\label{n}
\lim_{t\to\pm 1}\frac{f^{(n)}(t)}{(f'(t))^n}=0\quad\text{and}\quad
\lim_{t\to\pm 1}\frac{d}{dt}\frac{f^{(n)}(t)}{(f'(t))^n}=0.
\end{equation}

Consider the solid cylinder $D\times\Real$, where $D\subset\Real^2$
is the unit disc. On $D\times\Real$, consider the foliations with
the leaves $\mathbb S^1\times\Real$ and
$$L_\alpha=\{(x,f(|x|)+\alpha),\, x\in D\},\quad \alpha\in\Real,$$
where $|x|$ is the standard norm in $\Real^2$. We obtain the induced
foliation on the solid torus $D\times \mathbb S^1=D\times
(\Real/\mathbb Z)$. The projection $\Real\to\mathbb
S^1\subset\mathbb C$ is given by $\beta\mapsto e^{2\pi i\beta}$.

Consider the sphere $\mathbb S^3$ in $\Real^4=\Real^2\times\Real^2$
given by the equation $|z|^2+|w|^2=2$. Take two subsets $\mathbb
S^3_i$ $(i=1,2)$ of $\mathbb S^3$ defined by the conditions
$|z|^2\le |w|^2$ and $|z|^2\ge |w|^2$, respectively. It holds
$\mathbb S^3=\mathbb S^3_1\cup \mathbb S^3_2$. The map
$$\left(z,e^{2\pi i\beta}\right)\mapsto \left(z,\sqrt{2-|z|^2}e^{2\pi i \beta}\right)$$
defines a diffeomorphism of $D\times \mathbb S^1$ onto $\mathbb
S^3_1\subset \Real^2\times\mathbb C$. Take a function $h$ satisfying
the same properties as $f$; using it, construct a foliation on
$\mathbb S^1\times D$ and a diffeomorphism $\mathbb S^1\times D\to
\mathbb S^3_2$. In this way we obtain a Reeb foliation on $\mathbb
S^3$.

Consider the curve $$\gamma:(0,\sqrt{2})\to \mathbb
S^3\subset\Real^2\times\mathbb C,\quad x\mapsto \left(\frac{x}{\sqrt
2},\frac{x}{\sqrt 2},\sqrt{2-x^2}\right).$$ The curve $\gamma$ is
transverse to the Reeb foliation and it intersects at least once
each leaf of the foliation. For $x=1$, the curve passes the compact
leaf. For $0<x,y<1$, the points $\gamma(x)$ and $\gamma(y)$ belong
to the same leaf if and only if $f(x)=f(y)-n$ for some integer $n$.
Let
$$\varphi(x)=\left\{\begin{array}{ll}f^{-1}(f(x)+1),& \text{\rm if}\,\, 0<x<1;\\x, & \text{\rm if}\,\, x\geq 1.\end{array}
\right.$$ The function $\varphi$ is smooth, and it represents an
element of the holonomy group of the compact leaf. We see that the
leaf space $\mathbb{S}^3/\F$ may be identified with the orbit space
$$(0,\sqrt{2})/<\varphi,\psi>,$$ where the local diffeomorphism
$\psi $ is defined by the function $h$; here $<\varphi,\psi>$ is the
pseudogroup generated by the local diffeomorphisms $\varphi$ and
$\psi$.

Suppose that the first CL class of the foliation is trivial in
$H^2(S'(\mathbb{S}^3/\F))$. Then, on $S'((0,\sqrt{2}))$ there exists
a 1-form $\omega$ invariant under $\tilde\varphi$ and $\tilde\psi$
and such that
$$d\omega=dx_2\wedge dx_0.$$ Consider the new coordinate $t=f(x)\in(0,+\infty)$,
where $0<x<1$. Let $$\lambda=(\tilde
f^{-1})^*\omega\in\Omega^1(S'((0,+\infty))).$$ Since $\omega$ is
preserved by $\varphi$, $\lambda$ is preserved by the shift
$t_0\mapsto t_0+1$, i.e., $\lambda$ may be considered as an element
of $\Omega^1(S'(\mathbb S^1))$. It holds $$d(\lambda-t_2dt_0)=0.$$
Since the projection $S'(\mathbb S^1)\to \mathbb S^1$ is a homotopy
equivalence, there exist a constant $c$ and a function  $H$ on
$S'(\mathbb S^1)$ such that
$$\lambda-t_2dt_0=cdt_0+dH.$$
We obtain that for $x_0<1$ it holds
$$\omega=\tilde
f^*\lambda=\left(x_2+\frac{f''(x_0)}{f'(x_0)}\right)dx_0+cdf(x_0)+d\tilde
f^*H=x_2dx_0+d\left(\ln f'(x_0)+cf(x_0)+\tilde f^*H\right).$$ Since
the form $\omega$ is smooth for $x_0=1$, the function
\begin{equation}\label{smoothA}A(x_0,x_2,\dots)=\ln
f'(x_0)+cf(x_0)+H\left(f(x_0),\frac{x_2}{f'(x_0)}+\frac{f''(x_0)}{(f'(x_0))^2},\dots\right)\end{equation}
is smooth as well. Recall that the function $H$ is periodic in the
variable $t_0$. Let $x_2=x_3=\cdots=0$. Consider the sequence
$x_{0,n}=f^{-1}(n)$ and substitute it to \eqref{smoothA}. It holds
$x_{0,n}\rightarrow 1$ as $n\rightarrow +\infty$. Dividing the
equality \eqref{smoothA} by $f(x_{0,n})$, considering the limit
$n\rightarrow +\infty$, and using \eqref{coordtrans2}, \eqref{n}, we
get $c=0$. Considering \eqref{smoothA} and the same limit again, we
see that the sequence $\ln f'(x_{0,n})$ must be bounded. This
contradicts the properties of the function~$f$. \qed

\vskip0.2cm

The following corollary follows from Theorem \ref{ThCLReeb} by Proposition \ref{proptrivial1}.

\begin{cor} The  Godbillon-Vey-Losik class with values in
$H^3(S''(\mathbb S^3/\F)/\mathbb{Z})$ is non-trivial for the Reeb
foliations. \end{cor}

The   just proved theorem and its corollary imply that the GVL class
with values in $H^3(S''(M/\F)/\mathbb{Z})$ and the first CL class do
not satisfy the classical Mizutani-Morita-Tsuboi Theorem
\cite{Ref82} about triviality of the Godbillon-Vey class of
foliations almost without holonomy. Thus the  classes we consider
detect the non-trivial
 holonomy of the compact leaf. On the other hand, these classes do not contain the complete information
about the holonomy. This is shown by the following example.

\begin{ex} The Duminy-Losik class, the first Chern-Losik class and the Godbillon-Losik class are trivial for 1-dimensional orbifolds.
Indeed, consider the orbifold $[0,1)=(-1,1)/\mathbb{Z}_2$. Consider
the orbifold chart $$U=(-1,1)\to [0,1),\quad y\mapsto |y|.$$
Consider the coordinates \eqref{coord2}. The form
$$y_2dy_0=-y_1y_2 \theta_0$$ is invariant with respect to the transformation induced by the holonomy homomorphism $y\mapsto  -y$, and
the differential of this form equals $dy_2\wedge dy_0$, i.e., the
DL class and consequently the other classes are trivial. The
orbifold $[0,1]$ may be considered in the same way. This statement
follows also from Theorem \ref{ThDLclass}, since orbifolds admit
Riemannian metrics and consequently affine connections.
\end{ex}

Note that we may consider the Reeb foliations on the 2-dimensional
torus obtained by gluing two rings with 1-dimensional foliations. Each
such  foliation has the same holonomy group as a Reeb foliation on
the 3-dimensional sphere, consequently, the GVL class with values in
$H^3(S''(M/\F)/\mathbb{Z})$ and  the first CL class of this
foliation are non-trivial. On the other hand, the Godbillon-Vey
class of a foliation on a 2-dimensional manifold is always trivial.

Let $\F$ be a codimension-one foliation on a smooth manifold $M$.
Recall that a leaf $L$ from $\F$ is called resilient if for some
point $x\in L$ there exists a transversal $U$, a local diffeomorphism
$f$ defined on $U$ which represents an element of the holonomy group
at $x$, and another point $y\in L\cap U$ such that
$f^n(y)\rightarrow x$ as $n\rightarrow +\infty$. A theorem by Duminy
and Sergiescu \cite{Ref28} states that if the foliation $\F$ has no
resilient leaves, then its Godbillon-Vey class is zero. The case of
the Reeb foliations shows that the statement of this theorem does not
hold for the first CL class and for the GVL class with values in
$H^3(S''(M/\F)/\mathbb{Z})$.

It is clear how to define a resilient point of  a $\D_1$-space.
Just below, we give an example of a $\D_1$-space with a resilient point
and trivial characteristic classes.

\begin{ex} Let $X=\Real/<\varphi,\psi>$, where $\varphi(x)=x+1$ and
$\psi(x)=\frac{1}{2}x$. Then the orbit of the origin $0\in\Real$ is
the resilient point of $X$. The forms $x_2dx_0$ and $x_2 dx_0\wedge
dx_1$ are well-defined over $X$, consequently all characteristic
classes under consideration are trivial. \end{ex}

Recall that foliations with non-exponential growth have no resilient
leaves, see, e.g., \cite{Hurder00,Hurder16}.  Moussu and Pelletier,
and independently Sullivan conjectured that foliations with
non-exponential growth have trivial Godbillon-Vey class. A stronger
variant of this conjecture is proved by Hurder in \cite{Hurder86}.
We see that this conjecture is not true for the classes that we
consider.

Note also that each Reeb foliation on the 3-dimensional sphere is cobordant to zero \cite{Serg}, consequently, unlike the classical Godbillon-Vey class,  the classes under consideration are not invariant by cobordism.

Thus the (non-)triviality conditions for the classes under
consideration are quite different from the (non-)triviality
conditions for the Godbillon-Vey class and they are worth  studying.

\section{Sufficient condition for non-triviality}

Everywhere in the present section we will identify a vector field
$w\frac{d}{dx}$ on $\mathbb{R}$ with the function $w$. We denote
positive and non-negative real numbers by $\mathbb{R}_{>0}$ and
$\mathbb{R}_{\geq 0}$, respectively.

Let us recall the definition of the holonomy group for a $\D_n$-space $X$ at a point $p\in X$ \cite{L94}.  Let  $l:U\to X$, $k:V\to X$ be two $\D_n$-charts such that $l(u)=k(v)=p$, where $u\in U$ and $v\in V$. By the definition of a $\D_n$-space, there exists an open neighborhood $W$ of $u$ contained in $U$ and a morphism of $\D_n$-charts $f:W\to V$ such that $f(u)=v$. Denote by $G_{p,u,v}$ the set of germs at $u$ of such morphisms $f$. Let    $G_{p,u}=G_{p,u,u}$. For any other choice of the chart $l':U'\to X$ and point $u'\in U'$, the groups  $G_{p,u}$ and $G_{p,u'}$  are isomorphic. The group $G_{p,u}$ is called the holonomy group of 
$X$ at $p$. If $X=M/\F$, then this definition is equivalent to the usual definition of the leaf holonomy.

Let $X$ be a $\D_1$-space. Let $p \in X$ be a point with the
holonomy group $G$. If $G$ is non-trivial, then there exists a
 chart $k:V\to X$, an open subset $U\subset V$, a point $u\in U$ such that $k(u)=p$,   and a non-identity  morphism of charts $\varphi: U \rightarrow
V$ fixing $u$. We may assume that $\varphi$ is a diffeomorphism form $U$ onto its image. We denote the $n$-th
iteration of $\varphi$ by $\varphi_n = \varphi \circ \ldots \circ
\varphi$ ($\varphi_n$ is defined in an open neighborhood of the point $u$). We extend the last formula to all
 integers $n$ assuming that $\varphi_{-1} = \varphi^{-1}$ (for $n\geq 1$, the domain of definition of $\varphi_{-n}$ coincides
 with the image of $\varphi_n$). The following lemma is almost evident, but for completeness of the exposition we give the proof of it.

\begin{lem}
 Let $X$ be a $\D_1$-space and let $G$ be the holonomy group of $X$ at a point $p\in X$.  Then one of the following conditions holds:

i) $G = 1$;

ii) $G = \mathbb{Z}_2$;

iii) $G$ is infinite and contains the subgroup $\mathbb{Z}$
generated by some germ at $p$.

\end{lem}

{\bf Proof}. Assume that $G$ is non-trivial. If $G$ is finite, then we may consider  its elements as germs of local diffeomorphisms defined on a domain of a chart $k:U\to X$. Fix a point $u\in U$ such that $k(u)=p$. Consider the standard Riemannian metric $dx^2$ on $U$ and put
$$
d\bar{x}^2 = \frac{1}{|G|}\sum_{\gamma \in G}
(\varphi_\gamma^\prime)^2 dx^2,
$$
where $\varphi_\gamma$ is a local diffeomorphism representing $\gamma$.
The obtained metric is defined in a neighborhood of the point $u$ 
and is invariant with respect to $G$ in this neighborhood. This
immediately implies that $G = \O1 = \mathbb{Z}_2$ and $G$ is
generated by the germ of the local diffeomorphism $\varphi(x) = -x$.

Suppose that $G$ is infinite. Suppose that each element of $G$ is of finite order. By the above argument, each finite subgroup of $G$ is isomorphic to $\mathbb{Z}_2$. This implies that each non-identical element of $G$ is of order two. Consequently the group $G$ is Abelian. Finally, any two arbitrary non-identical elements generate a finite subgroup of $G$, which is isomorphic to $\mathbb{Z}_2$, i.e., these elements coincide. This implies $G=\mathbb{Z}_2$ and contradicts the infiniteness of $G$.
 Therefore,
there exists a morphism of infinite order, and its germ at $u$ generates the subgroup $\mathbb{Z}$ in~$G_{p,u}$.  \qed

Let $X$ be a $\D_1$-space and $p\in X$. Let $k:V\to X$ be a chart, $U\subset V$ an open subset, $u\in U$, and $k(u)=p$. Let $\varphi:U\to V$ be a morphism of charts.
The point $u$ is called a fixed non-hyperbolic point of $\varphi$ 
if $\varphi(u)=u$ and $|\varphi^\prime(u)| = 1$ \cite{L90}. 
In this case we say that the element of the holonomy group at the point $p$ corresponding to $\varphi$ is non-hyperbolic.
We call the point  $p \in X$ non-hyperbolic
if the holonomy group at the point $p$ contains an \emph{\bf infinite order} non-hyperbolic element.

We formulate now the main theorem of this section.

\begin{theorem} \label{Nontrivil1A}
    Suppose that a $\D_1$-space $X$ has a non-hyperbolic point $p$. Then the first Chern-Losik class of $X$ is
non-trivial in $H^2(S'(X))$.
\end{theorem}

\begin{cor}\label{cornontrivil}
 Suppose that a $\D_1$-space $X$ has a non-hyperbolic point. Then the  Godbillon-Vey-Losik class is non-trivial in
$H^3(S''(X)/\mathbb{Z})$.
\end{cor}

The next example shows that the statement of Theorem
\ref{Nontrivil1A} is generally not true for spaces with hyperbolic
points.

\begin{ex}
Suppose that a $\D_1$-space $X$ is the orbit space
$\mathbb{R}/\langle \varphi \rangle$, where $\varphi$ is a
 diffeomorphism with a unique fixed point which is hyperbolic. Then the first Chern-Losik class of $X$ is
trivial.

Indeed,  from the results of \cite{Sternberg} it follows that by a
smooth change of the coordinate,  the diffeomorhism $\varphi$ can be
conjugated to a linear diffeomorphism, i.e.,
$$
\varphi(x) = kx,
$$
where $k \neq 0, \pm 1$. Then the $1$-form $\omega = x_2 dx_0$ is
invariant with respect to the  lift $\tilde\varphi:S'(X)\to S'(X)$ of $\varphi$,  and $d \omega$ represents the trivial
class $\CL_1(X)$.

Remark that Losik \cite{L90} proved that if $\varphi$ is a
diffeomorphism of $\Real$ such that $\varphi(x+1)=\varphi(x)+1$ and
it has two fixed points with different values of $\varphi'$ at these
points (i.e., it is possible that the both points are hyperbolic), then the cohomology class
of the form $\theta_2\wedge \theta_0$ is non-trivial in
$H^2(S'(\Real/(\mathbb{Z}\times <\varphi>)))$.
\end{ex}

Let $X$ be as in the statement of Theorem~\ref{Nontrivil1A}. Consider a chart $k:V\to X$, an open subset $U\subset V$, and a point $u\in U$ such that $k(u)=p$. Let $\varphi:U\to V$ be a morphism of charts representing the infinite order element of the holonomy group at $p$ with the fixed non-hyperbolic point $u$.
 Let us denote the set of all fixed points of  $\varphi$ by $F = F(\varphi)$.  We say that a point $z \in F$ is right (resp., left)
  semi-isolated if and only if the semi-interval $I = (z,z+\varepsilon)$ (resp.,  $I=(z-\varepsilon, z)$) does not contain  points from $F$ for
   sufficiently small $\varepsilon>0$. In this situation we will call the interval $I$ a right (resp., left)  open segment of $z$, and we say that the corresponding point in $X$ is semi-isolated.

First we prove a weak version of Theorem \ref{Nontrivil1A}.

\begin{theorem}\label{Nontrivil1}
Suppose that a $\D_1$-space $X$ has a semi-isolated non-hyperbolic point. Then the first Chern-Losik class
of $X$ is non-trivial in $H^2(S'(X))$.
\end{theorem}

{\bf Proof of Theorem \ref{Nontrivil1}.}  By the assumption, there exists a chart $k:V\to X$, an open subset $U\subset V$, and a morphism $\varphi:U\to V$ with a    semi-isolated non-hyperbolic point  $u\in U$ such that $\varphi$ represents an element of infinite order in the holonomy group at the point $k(u)\in X$.
 Without loss of generality we may assume
that  $U=V=\mathbb{R}$, $u=0$, $\varphi(0)=0$, and $\varphi(x)>x$ for
all $x>0$ (if $\varphi(x)<x$, then we replace $\varphi$ by
$\varphi^{-1}$). From results of \cite{Szekeres},
\cite{Kopell} it follows  that there exists a one-parameter group
$\varphi_t:\mathbb{R}_{\geq 0} \rightarrow \mathbb{R}_{\geq 0}$
extending the iterations $\varphi_n$ and generated by a vector field
$v$ on $\mathbb{R}_{\geq 0}$ of class $C^1$. Following \cite{Eyn},
we call $v$ the Szekeres field of $\varphi$. Let us fix some $\xi
>0$ and consider the new coordinate $t \in \mathbb{R}$ defined   by the relation $x =
\varphi_t(\xi)$. Denote this coordinate transformation by $x=g(t)$
and let $g=f^{-1}$. It is clear that $f:\mathbb{R}_{>0} \rightarrow
\mathbb{R}$ and $g: \mathbb{R} \rightarrow \mathbb{R}_{>0}$ are
diffeomorphisms of class $C^2$ and it holds
$$
g^\prime(t) = v(g(t)),\quad f^\prime (x) = \frac{1}{v(x)}.
$$

Let us prove some lemmas.

\begin{lem}\label{LemmOnVectorField}
Let $f$ have continuous derivatives up to order $n+1$ at a point
$x>0$. Then $v$ has continuous derivatives up to order $n$ at $x$
and it holds
\begin{equation}\label{EqLemmOnvectorField}
v^{(n)}(x)(v(x))^{n-1}= -\frac{f^{(n+1)}(x)}{(f^\prime(x))^{n+1}} +
Q_n\left(\frac{f^{\prime \prime}(x)}{(f^\prime(x))^{2}}, \ldots,
\frac{f^{(n)}(x)}{(f^\prime(x))^{n}}\right),\quad n\geq 1,
\end{equation}
where $Q_n(u_1, \ldots, u_{n-1})$  is a polynomial with integer
coefficients of degree less or equal to $n$, and $Q_1 = 0$.
\end{lem}

{\bf Proof}. For the proof we use induction over $n$. The first
step of the induction is trivial. Suppose that the formula is true
for some natural number $n$. Then, differentiating the equation
\eqref{EqLemmOnvectorField}, we obtain
$$
v^{(n+1)}(x)(v(x))^{n-1} + (n-1)v^{(n)}(x)(v(x))^{n-2}v^{\prime}(x)
=
$$
$$
-\frac{f^{(n+2)}(x)}{(f^\prime(x))^{n+1}} + (n+1)\frac{f^{(n+1)}(x)
f^{\prime \prime}(x)}{(f^\prime(x))^{n+2}} +
$$
$$
\sum_{k=1}^{n-1} \frac{\partial Q_n}{\partial u_k}
\left(\frac{f^{\prime \prime}(x)}{(f^\prime(x))^{2}}, \ldots,
\frac{f^{(n)}(x)}{(f^\prime(x))^{n}}\right) \left(
\frac{f^{(k+2)}(x)}{(f^\prime(x))^{k+1}} - (k+1)\frac{f^{(k+1)}(x)
f^{\prime \prime}(x)}{(f^\prime(x))^{k+2}}   \right).
$$
Multiplying now both sides of the last equality by
$v(x)=1/f^\prime(x)$, we get the formula we need with
$$
Q_{n+1}(u_1,\ldots, u_n) = 2 u_1 u_{n} + (n-1)u_1
Q_n(u_1,\ldots,u_{n-1}) +
$$
$$
\sum_{k=1}^{n-1}\frac{\partial Q_n}{\partial
u_k}(u_1,\ldots,u_{n-1})(u_{k+1} - (k+1)u_1 u_{k}).
$$
Finally, $\deg(Q_2)=2$ and $\deg(Q_{n+1}) \leq \deg(Q_n) + 1$ for $n
\geq 2$, which implies $\deg(Q_n) \leq n$ for all~$n$. \qed

\begin{lem}\label{vProperties}
The  Szekeres vector field has the following properties:

i) $v(0) = 0$;

ii) if $\varphi^\prime(0)=1$, then $v^\prime(0) = 0$.
\end{lem}

{\bf Proof}. Suppose that for some $t \in \mathbb{R}$ it holds
$\varphi_t(0)>0$. Then $\varphi_r(0)>0$ for some rational $r \in
\mathbb{Q}$. This implies that $\varphi_\frac{1}{n}(0)>0$ for some
$n \in \mathbb{N}$. We obtain that $\varphi(0) >0$, which gives a
contradiction. Consequently, $\varphi_t(0)=0$ for every $t \in
\mathbb{R}$. By definition of $v$ we have
\begin{equation}\label{EqV1}
\frac{d}{dt}\varphi_t(x) = v(\varphi_t(x)).
\end{equation}
Substituting $x=0$ to \eqref{EqV1}, we immediately obtain $v(0) =
0$.

Now, for every $n \in \mathbb{N}$, the equation
$$
\varphi^\prime(0) = (\varphi_\frac{1}{n}^\prime(0))^n = 1
$$
implies $\varphi_\frac{1}{n}^\prime(0) = 1$. Then $\varphi_{r}(0) =
1$ for all $r \in \mathbb{Q}$ and consequently $\varphi_t^\prime(0)
= 1$ for all $t$. Differentiating \eqref{EqV1} by $x$ and
substituting $x=0$, we obtain $v^\prime(0)=0$. \qed

\vskip0.3cm

The next lemma generalizes the  proof of  non-triviality of the
first CL class for the Reeb foliations given in the previous section.

\begin{lem}\label{ReebGeneralization}

Suppose that a $\D_1$-space $X$ has a right (resp., left)
semi-isolated point $p \in X$ with the holonomy group at $p$ containing the germ of a local diffeomorphism with the Szekeres vector
field $v$ satisfying the following properties: $v$ is
$C^\infty$-smooth in a right (resp., left)  open segment $I$ of $p$,
and for each $n\in \mathbb{N}$ there exists a constant $C_n$ such
that
\begin{equation}\label{ZeroCondition}
v^\prime(0) = 0,
\end{equation}
\begin{equation}\label{VectorFieldCondition}
\left|v^{(n)}(p)(v(p))^{n-1}\right| \leq C_n,\quad p \in I
\end{equation}
Then the first Chern-Losik class of $X$ is non-trivial.
\end{lem}

{\bf Proof}. As above we assume that $U =
\mathbb{R}$, $u=0$, $\varphi(x)<x$ for $x>0$. Consider the Szekeres
field $v$, the one-parameter group $\varphi_t$ and the new
coordinate $t$ as above. Assume that the first CL class is trivial.
Then there exists a $\tilde{\varphi}$-invariant $1$-form $\omega$
satisfying the property
$$
d \omega = dx_2 \wedge dx_0, 
$$ here $\tilde \varphi:S'(\Real)\to S'(\Real)$ is the lift of $\varphi$.
The equality
$$
\lambda = \tilde{g}(\omega)
$$
defines a $\widetilde{ g^{-1}\circ\varphi\circ g}$-invariant
$1$-form on $S'(\mathbb{R})$. Then the $1$-form $\lambda - t_2\wedge
dt_0$ is closed and the Poincar\'{e} Lemma implies that $$\lambda -
t_2\wedge dt_0 = dH$$ for some smooth function $H = H(t_0, t_2,
\ldots, t_q)$ on $S^\prime(\mathbb{R})$ ($q \in \mathbb{N}$). It is
easy to see that $$ (g^{-1}\circ \varphi \circ g) (t) =
f(\varphi(\varphi_t(\xi))) = f (\varphi_{t+1}(\xi)) = t+1$$ and
$$\widetilde{(g^{-1}\circ \varphi \circ g)} (t_0, t_2, \ldots, t_q) = (t_0 + 1, t_2,
\ldots, t_q).$$ Therefore, $$dH(t_0, t_2, \ldots, t_q )=dH(t_0+1,
t_2, \ldots, t_q)$$ and $$H(t_0+1, t_2, \ldots, t_q) = H(t_0, t_2,
\ldots, t_q ) + c_0$$ for some constant $c_0$. The condition
\eqref{VectorFieldCondition} and Lemma \ref{LemmOnVectorField} imply
that there exists a constant $c_1$ such that
\begin{equation}\label{Ineq1}
\left|\frac{f^{(n)}(p)}{(f^\prime(p))^n}\right| \leq c_1
\end{equation}
for all $n \leq q$ and all $p \in I$. Since the function $H$ is
smooth, there exists a constant $c_2$ such that for every $0 \leq
t_0 \leq 1$ and for all $t_2, \ldots, t_q$ satisfying
$$
|t_n|\leq c_1,\quad 2 \leq n \leq q
$$
it holds
\begin{equation}\label{Ineq2}
|H(t_0, t_2, \ldots, t_q)| \leq c_2.
\end{equation}
Returning to the coordinate $x$, we obtain
$$
\omega = \tilde{f} (\lambda) = x_2 \wedge dx_0 + d \left( \ln
|f^\prime(x_0)| + H\left( f(x_0), \tilde{f}(x_2), \ldots, \tilde{f}
(x_q) \right) \right).
$$
Since $\omega$ is smooth on $S^\prime (\mathbb{R}_{\geq 0})$, the
function
$$
A(x_0, x_2, \ldots, x_q) = \ln |f^\prime(x_0)| + H\left( f(x_0),
\tilde{f}(x_2), \ldots, \tilde{f} (x_q) \right)
$$
is smooth on $S^\prime (\mathbb{R}_{\geq 0})$. Then the function
$$
A(x_0) = A(x_0, 0, \ldots, 0) = \ln |f^\prime(x_0)| + H\left(
f(x_0), \frac{f^{\prime \prime}(x_0)}{(f^\prime(x_0))^2}, \ldots,
\frac{f^{(q)}(x_0)}{(f^\prime(x_0))^q} \right)
$$
is bounded on $I$. We have
\begin{equation}\label{equ1}
A(x_0) = \ln |f^\prime(x_0)| + H\left( f(x_0) - [f(x_0)],
\frac{f^{\prime \prime}(x_0)}{(f^\prime(x_0))^2}, \ldots,
\frac{f^{(q)}(x_0)}{(f^\prime(x_0))^q} \right) + c_0 [f(x_0)]
\end{equation}
(here we denote by $[y]$ the integer part of $y$). It follows from
(\ref{Ineq1}) and (\ref{Ineq2}) that the second summand on the right
hand side of (\ref{equ1}) is bounded by $c_2$ on $I$. We have
$f(x_0) \rightarrow \infty$ and $f^\prime(x_0) = 1/v(x_0)
\rightarrow -\infty$ as $x_0 \rightarrow 0$. Taking into the account
the facts  that the functions $\ln f^\prime (p)$ and $f(p)$ are
differentiable for $p>0$, and using \eqref{ZeroCondition}, we obtain
$$
\lim_{x_0\rightarrow 0} \frac{\ln f^\prime(x_0)}{f(x_0)} =
\lim_{x_0\rightarrow 0} \frac{f^{\prime
\prime}(x_0)}{(f^\prime(x_0))^2} = -v^\prime(0) = 0.
$$
Dividing \eqref{equ1} by $f(x_0)$, taking the limit $x_0 \rightarrow
0$ and using the last equality, we immediately obtain that $c_0 =
0$. Now, taking the upper limit $x_0 \rightarrow 0$ in \eqref{equ1}
again, we obtain that $\ln |f^\prime(x_0)|$ is bounded on $I$, which
 gives a contradiction. \qed

We complete now the proof of Theorem \ref{Nontrivil1}. First
consider the case of a  diffeomorphism $\varphi$ which is not tangent
to the identity. It follows from \cite{Takens} that the vector field
$v$ is of class $C^\infty$ on $\mathbb{R}_{\geq 0}$. Then $v^{(n)}
(v^\prime)^{n-1}$ is smooth for all $n\geq 1$ and the condition
\eqref{VectorFieldCondition} holds true. The condition
\eqref{ZeroCondition} follows from Lemma \ref{vProperties}. Now
Lemma \ref{ReebGeneralization} completes the proof of the theorem in
this case.

The remaining case of the diffeomorhisms infinitely tangent to the
identity is  the most interesting. For that case, in \cite{Serg} it
is proved that $v$ is of class $C^\infty$ on $\mathbb{R}_{>0}$.
Consider set $T$ consisting of points $t \in \mathbb{R}$ such that
$\varphi_t$ is $C^\infty$-smooth. It is clear that $\mathbb{Z}
\subset T \subset \mathbb{R}$. In \cite{Serg}, a diffeomorhism
$\varphi$ is constructed  with $T=\mathbb{Z}$ (see also
\cite{EynardETD}). In \cite{Eyn}, it is shown that $T$ is a Cantor
set between $\mathbb{Z}$ and $\mathbb{R}$. In the both cases, the
Szekeres field $v$ is not of class $C^2$ on $\mathbb{R}_{\geq 0}$.
Nevertheless, one can show that the function $v^{(n)}v^{n-1}$ is
bounded on $\mathbb{R}_{\geq 0}$.

To prove it, we need some further results from \cite{Serg}. Following
\cite{Serg}, let
$$
\Delta x = x - \varphi(x),\quad \Delta_0 x = \sup_{0 \leq y \leq x}
\Delta y.
$$
In \cite[Lemma 3.6]{Serg} it is proved  that
$$
|v^{(n)}(x)(v(x))^{n-1}| = O\left((\Delta_0 x)^{n - \eta}\right))
$$
for every real $\eta >0$. We have $\Delta_0 x = o(1)$, so
$|v^{(n)}v^{n-1}|$ is bounded in a neighbourhood of $x=0$. Now the
argument from the case of a diffeomorphism not tangent to the
identity may be applied, and this completes the proof of Theorem
\ref{Nontrivil1}. \qed

{\bf Proof of Theorem \ref{Nontrivil1A}.} By the assumption, there
exists a chart $k:V\to X$, an open subset $U\subset V$, a point $u\in U$ such that $k(u)=p$, and a morphism $\varphi:U \rightarrow V$ of charts that satisfies $\varphi(u)=u$, $|\varphi^\prime(u)| = 1$ and represents 
an infinite order element in the holonomy group $G$ at the point
$p$. Let $$S = \{ z \in U |
\varphi(z) = z, |\varphi^\prime(z)| = 1 \}.$$ The subset $S\subset U$ is
closed, non-empty and proper in $U$ (if $S=U$, then $|\varphi^\prime|
= 1$ on $U$ and the germ of $\varphi$ at $u$ has finite order in $G_{p,u}$). Then $U
\backslash S$ is a non-empty union of open intervals. This implies
that there exists a $z$ in $S$ such that $\varphi(z)=z$, $|\varphi'(z)|=1$, and the point $z$ is  semi-isolated. By the definition of the holonomy group, $\varphi$ represents an element of infinite order in the holonomy group at the point $k(z)\in X$. By theorem \ref{Nontrivil1}, the first CL class of $X$ is non-trivial.  \qed

\section{Codimension-one foliations without holonomy}

 Let
$M$ be a compact connected smooth manifold and let $\F$ be a
codimension-one foliation without holonomy on $M$.
 Novikov \cite{Novikov} proved that the universal covering $\tilde{M}$ of $M$ is diffeomorphic to
 $\tilde{L}\times \mathbb{R}$, where $\tilde{L}$ is the universal covering of any leaf of $F$. Moreover, the deck transformations of
  $\tilde{M}$ preserve the leaves $\tilde{L}\times \{z\}$ and induce the following homomorphism with an Abelian
  image:
$$
q: \pi_1(M) \rightarrow \Diff_{+}({\mathbb R}),
$$
where $\Diff_{+}({\mathbb R})$ is the group of orientation
preserving $C^\infty$-diffeomorphisms of $\mathbb{R}$. More
precisely, for every $\alpha \in \pi_1(M)$ and $x \in \tilde{M}$,
the homomorphism $q$ is defined by the relation $$\pi(\alpha(x)) =
q(\alpha)p(x),$$ where $p: \tilde{M}=\tilde{L} \times {\mathbb R}
\rightarrow {\mathbb R}$ is the projection onto the second factor,
and $\alpha(x)$ is the deck transformation applied to $x$. It is
important to note that since the foliation $\F$ is without
holonomy, the non-trivial elements of $\Im(q)$ have no fixed
points.

Consider the free action of the group $\pi_1(M)$ on $\tilde{M}
\times {\mathbb R}$ as following:
$$
\alpha \cdot (x,z) = (\alpha(x), q(\alpha)z),
$$
where $\alpha \in \pi_1(M)$, $x \in \tilde{M}$, $z \in {\mathbb R}$.
The quotient manifold $E=\pi_1(M) \backslash \tilde{M} \times
{\mathbb R}$ has a natural foliation defined by the leaves
$\tilde{M} \times \{z\}$. Let us denote this foliation by $\tilde
\F$. Thus for the leaf space of this foliation we get
$$E/\tilde\F=\Real/\Im(q).$$ In \cite{Ref84}, a
cross-section $\sigma: M \rightarrow E$ transversal to the leaves of
$\tilde \F$ is defined in such a way that the foliation $M$ induced
by $\sigma$ coincides with $\F$. Until now we have been following the
ideas of Morita and Tsuboi from \cite{Ref84}.

Now, the section $\sigma$ induces the morphism of the leaf spaces
$$M/\F\to E/\tilde\F.$$
From the existence of this map and the naturality of the
characteristic classes it follows that if some class is trivial
for $E/\tilde\F$, then this class is trivial for~$M/\F$ as well.

The factorization map $$k:U=\Real\to \Real/\Im(q)=E/\tilde\F$$
provides a $\D_1$-chart on $E/\tilde\F$. The elements of $\Im(q)$
are morphisms of this chart. In this way we get a $\D_1$-atlas on
the space $E/\tilde\F$.

Consider pairwise commuting generators $\varphi_0, \varphi_1,
\ldots, \varphi_p$ of $\Im(q)$. It is easy to see that we can choose
the coordinate $z$ on $\mathbb{R}$ in such a way that $\varphi_0$
acts as the unit shift: $\varphi_0(z) = z+1$. Then the
$\varphi_0$-action on $\mathbb{R}$ generates the circle
$\mathbb{S}^1=\mathbb{R}/\mathbb{Z}$. Since the diffeomorphisms $\varphi_i$,
$i=1,\ldots,p,$ commute with $\varphi_0$, they may be considered
 as
$C^\infty$-diffeomorphisms of $\mathbb{S}^1$. Note that the obtained
diffeomorphisms of $\mathbb{S}^1$ are again pairwise commuting and have no
fixed points. Even more, by the assumption, the non-trivial elements of the group generated by $\varphi_1,
\ldots, \varphi_p$ have no fixed points. In particular, non-zero powers of each element $\varphi_i$ $i=1,\ldots,p,$ have no fixed points. This implies that $\varphi_1,
\ldots, \varphi_p$ have no periodic points on~$\mathbb{S}^1$. 

Consider a $C^\infty$-diffeomorphism $\varphi: \mathbb{S}^1 \rightarrow \mathbb{S}^1$
without fixed points. Recall that the rotation number of $\varphi$
is equal to the following limit:
$$
\alpha = \rho(\varphi) = \lim_{n\rightarrow \infty}
\frac{\varphi^n(z)}{n} \ \mod \  1.
$$
In what follows, we change the moduli reduction in this definition
assuming that $\alpha \in (0,1]$. It is a well-known fact that the
limit exists and does not depend on the choice of $z \in \mathbb{S}^1$,  and the
absence of periodic points implies the irrationality of $\alpha$, see, e.g.,~\cite{Navas}. According to
the classical Theorems by Poincar\'e and Denjoy \cite{Navas}, the diffeomorphism $\varphi$ is
topologically conjugate to the rotation by the angle $\alpha$, that
is, there exists a homeomorphism $f:\mathbb{S}^1 \rightarrow \mathbb{S}^1$ such that
$$f(\varphi(z)) = f(z)+ \alpha \ \mod 1.$$ If $f$ is  of class
$C^\infty$, then $\varphi$ is said to be $C^\infty$-conjugate to
a rotation. It is easy to see that because of the pairwise
commutativity of $\varphi_1, \ldots, \varphi_p$, all
diffeomorphisms $\varphi_i, i=1, \ldots, p,$ are
$C^\infty$-conjugate to rotations if and only if this is true for
a particular $\varphi_i$. Moreover, let $\rho(\varphi_i,
\varphi_j)$ be the rotation number of $\varphi_j$ with respect to
the circle generated by $\varphi_i$. Then it is obvious that
$$\rho(\varphi_i, \varphi_j) = \frac{1}{\rho(\varphi_j,
\varphi_i)},$$ and the property to be $C^\infty$-conjugate does
not depend on the choice of the diffeomorphism $\varphi_i$
generating the circle. This shows that the following definition is
correct. We say that $\Im(q)$ is $C^\infty$-conjugate to shifts
if $\Im(q)$ is generated by diffeomorphisms $\varphi_0, \ldots,
\varphi_p$, and $\varphi_1$ is $C^\infty$-conjugate to a rotation
of circle generated by $\varphi_0$.

The goal of this section is to prove the following theorem.

\begin{theorem}\label{Thwouthol} Let $M$ be a compact connected manifold with a
codimension-one foliation $\F$ without holonomy and let $q: \pi_1(M)
\rightarrow \Diff_+(\mathbb{R})$ be the homomorphism defined above.
If $\Im(q)$ is $C^\infty$-conjugate to shifts, then the
Duminy-Losik, the first Chern-Losik class, and the
Godbillon-Vey-Losik class of the foliation $\F$ are trivial.
\end{theorem}

Before we prove the theorem let us discuss the sufficient
condition from the statement of the theorem in more details. A
number $\alpha \in \mathbb{R}$ is said to satisfy the Diophantine
condition if there exists $\beta\geq 0$ provided that, for some
constant $C=C(\alpha,\beta)>0$, the estimate
$$
\left| \alpha - \frac{p}{q} \right| > \frac{C}{q^{2+\beta}}
$$
holds for every rational fraction $p/q \in \mathbb{Q}$. In this
situation the number $\beta$ is called the exponent of $\alpha$. The
following fact was proved in \cite{Yoc}, see also 
\cite{Katok,Katznelson}. Suppose that $\varphi$ is a $C^r$-diffeomorphism of
the circle with $r$ real and greater than $1$. Suppose further that
$\varphi$ has degree $1$ and its rotation number satisfies the
Diophantine condition with an exponent $\beta$. If $r>2+\beta$ and
$\varepsilon>0$, then the conjugacy between the rotation by
$\rho(\varphi)$ and $\varphi$ is of class $C^{r-1-\beta
-\varepsilon}$.

\begin{cor} In the notation of Theorem 1, if
$\rho(\varphi_i, \varphi_j)$ is Diophantine for some $i\neq j$, then
the Duminy-Losik class,
 the first Chern-Losik class, and the Godbillon-Vey-Losik class of
 the
foliation $\F$ are trivial.\end{cor}

Non-Diophantine numbers are called Liouville numbers and have zero
Hausdorff dimension \cite{Oxtoby}. This shows that the case of all
 numbers $\rho(\varphi_i, \varphi_j)$, $i\neq j$, being non-Diophantine,
is rather peculiar.

{\bf Proof of Theorem \ref{Thwouthol}}. We need to prove 
triviality of the characteristic classes for the $\D_1$-space
$E/\tilde \F=\Real/\Im(q)$. Assuming that
$\mathbb{S}^1=\mathbb{R}/<\varphi_0>$, we can find an orientation preserving
$C^{\infty}$-diffeomorphism $f: \mathbb{S}^1 \rightarrow \mathbb{S}^1$ such that for
every $z \in \mathbb{S}^1$, $i=1,\ldots, p$, it holds
$$
f(\varphi_i (z)) = f(z) + \alpha_i \ \mod  1,
$$
where $\alpha_i = \rho(\varphi_i)$. We can lift $f$ to a
diffeomorphism $F:\mathbb{R} \rightarrow \mathbb{R}$ with the
property $F(z+1) = F(z) + 1$. Using the $1$-periodicity of
$F^\prime$, define the $C^{\infty}$-function $\lambda:\mathbb{S}^1
\rightarrow \mathbb{R}$ as  follows:
$$
\lambda(z) = - \ln F^\prime(z) = - \ln f^\prime(z),\  z \in \mathbb{S}^1.
$$
Then
$$
(\varphi_i^*\lambda)(z) = - \ln f^\prime(\varphi_i(z)) = - \ln
(f\circ\varphi_i)^\prime(z) + \ln \varphi_i^\prime(z) = \lambda(z) +
\ln \varphi_i^\prime(z).
$$
This implies that the function $\frac{d\lambda}{dz}$ defines an
affine connection on $E/\tilde \F=\Real/\Im(q)$. Consequently, the
Duminy-Losik class  of $E/\tilde \F$ is trivial. This implies the
proof of the theorem.  $\Box$

\section{Conclusion} We have considered above various
characteristic classes for a codimension-one foliation $\F$ on a
smooth manifold $M$. In \cite{CM,Gal17}, one may find additional
 information about the Godbillon-Vey class with values in the
\v{C}ech-de~Rham cohomology $\check H^3(M/\F)$ of the leaf space, this cohomology is isomorphic to the cohomology of the classifying space of the holonomy groupoid of the foliation.
The following diagram follows from the results obtained above and in \cite{CM,BGG}. The diagram illustrates the triviality conditions for all
these classes

$$
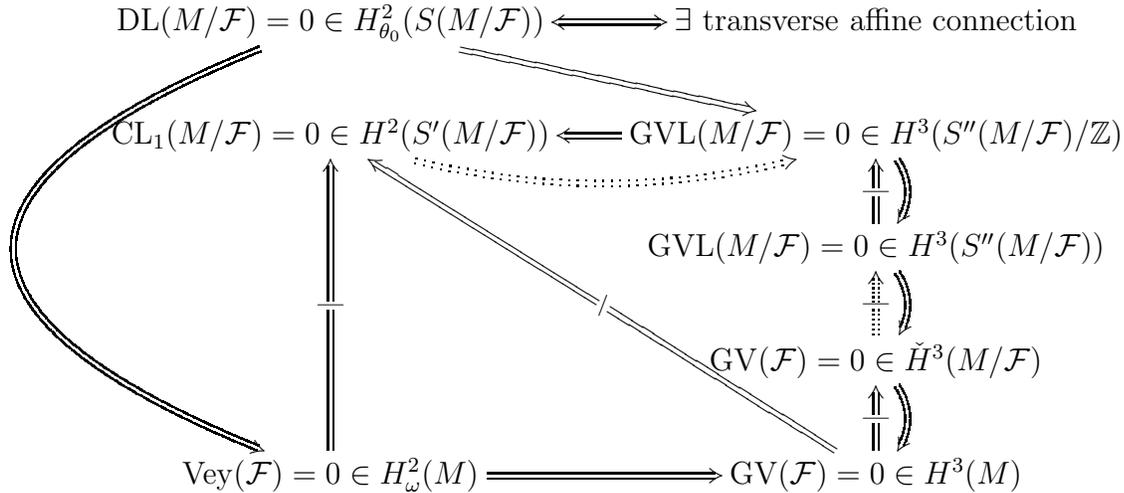
\xymatrix{
{\rm DL}(M/\F)=0\in H^2_{\theta_0}(S(M/\F))\ar@{=>}[dr]
\ar@/_10pc/@{=>}[dddd]  
\ar@{<=>}[r] & \exists\,\,\text{\rm transverse affine
connection}\\
{\rm CL}_1(M/\F)=0\in H^2(S'(M/\F)) \ar@/_1.5pc/@{:>}[r]\ar@{<=}[r] & {\rm GVL}(M/\F)=0\in H^3(S''(M/\F)/\mathbb{Z}) \\
&  {\rm GVL}(M/\F)=0\in H^3(S''(M/\F))\ar@/_1pc/@{<=}[u] \ar@{=>}[u]|{\text{\bf ---}}\\
&  {\rm GV}(\F)=0\in\check H^3(M/\F) \ar@/_1pc/@{<=}[u] \ar@{:>}[u]|{\text{\bf ---}}\\
{\rm Vey}(\F)=0\in H^2_\omega(M)  \ar@{=>}[uuu]|{\text{\bf ---}} \ar@{=>}[r] & {\rm GV}(\F)=0\in
H^3(M)\ar@{=>}[uuul]|{/} \ar@/_1pc/@{<=}[u] \ar@{=>}[u]|{\text{\bf ---}} \\
}$$

 The same diagram is valid if we consider the second order frame bundles, i.e., if we change all $S$ by $S_2$; 
in that case  the dot implications $:::>$ are valid  .

Thus the triviality of the Duminy-Losik class is the strongest
condition, which is equivalent to the existence of a transverse
affine connection. The triviality of the classical Godbillon-Vey
class is the weakest condition, and it was studied by many authors,
see the reviews \cite{Hurder00,Hurder16} or the book \cite{FolII}.
In particular, the Godbillon-Vey class of a foliation without resilient leaves is trivial, and information about  geometry of such foliations should give some other classes. 

For the Reeb foliations, the Vey class, the classical Godbillon-Vey class, and the Godbillon-Vey class with values in $\check H^3(M/\F)$ are trivial (the triviality of the last class follows directly from results in non-commutative geometry \cite{Connes,H-k}, see also \cite{Gal17}). The first Chern-Losik class and the Godbillon-Vey-Losik class with values in $H^3(S''(M/\F)/\mathbb{Z})$ are non-trivial for the Reeb foliations, and these classes detect the compact leaf with non-trivial holonomy. The Godbillon-Vey-Losik class with values in $H^3(S''(M/\F))$ turns out to be more delicate: it is trivial for some Reeb foliations and non-trivial for some other Reeb foliations. By that reason it should contain rich information about geometry and dynamics of the foliation leaves.

{\bf Acknowledgements.} The authors are thankful to Steven Hurder
for useful email communications and to the anonymous referee for valuable comments and suggestions that have significantly
improved the manuscript. The work was supported by grant no.
18-00496S of the Czech Science Foundation. Ya. V. Bazaikin was partially supported by the Program of Fundamental Scientific Research of the SB RAS No. I.1.2., Project No. 0314-2019-0006


\end{document}